\documentclass[10pt, a4paper]{article}

\usepackage{amsthm, amssymb, amsmath, textcomp, mathrsfs}
\usepackage[latin1,ansinew]{inputenc}
\usepackage[pdftex]{graphics}
\usepackage{graphicx}
\usepackage{longtable}
\usepackage{scrpage2}
\usepackage{hyphenat}
\usepackage{paralist}
\usepackage{color}
\usepackage{extpfeil}

\makeatletter
\renewcommand{\section}{%
  \@startsection {section}{1}{\z@}%
                 {-3.5ex plus -1ex minus -.2ex}%
                 {2.3ex plus.2ex}%
                 {\normalfont\Large\bfseries}}

\renewcommand{\subsection}{%
  \@startsection {subsection}{1}{\z@}%
                 {-3.5ex plus -1ex minus -.2ex}%
                 {2.3ex plus.2ex}%
                 {\normalfont\large\bfseries}}
                 
\renewcommand{\subsubsection}{%
  \@startsection {subsubsection}{1}{\z@}%
                 {-2.5ex plus -1ex minus -.2ex}%
                 {1.3ex plus.2ex}%
                 {\normalfont\bfseries}}
\makeatother

\setheadsepline{.4pt}
\clearscrheadfoot
\lehead{\parbox[b]{0.8cm}{\pagemark}\textsc{Chapter \thechapter}}
\rohead{\textsc{\rightmark}\parbox[b]{0.8cm}{\hfill \pagemark}}

\newcommand{\Completion}{May 22, 2012}

	\allowdisplaybreaks
	
	\newcommand{\dStr}[1]{{\rm \textbf{#1}}}	

	\usepackage{stmaryrd}

	\usepackage{hyperref}
	\hyperbaseurl{.}	

	\newtheorem{definition}{Definition}
	\newtheorem{theorem}{Theorem}
	\newtheorem*{unnumberedtheorem}{Theorem}
	\newtheorem{remark}{Remark}
	\newtheorem{corollary}{Corollary}
	\newtheorem{lemma}{Lemma}
	\newtheorem{example}{Example}	

	\newcommand{\grad}{\operatorname{grad}}
	
	\newcommand{\spann}{\operatorname{span}}

	\newcommand{\trace}{\operatorname{trace}}

	\newcommand{\pr}{\operatorname{pr}}	
	\newcommand{\To}{\longrightarrow}

	\newcommand{\Hol}{\operatorname{Hol}}
	
	\newcommand{\Isom}{\operatorname{Isom}}
		
	\newcommand{\LAkill}{\mathfrak{kill}}
	\newcommand{\he}{\mathfrak{he}}
	\newcommand{\He}{\operatorname{He}}
	\newcommand{\End}{\operatorname{End}}

	\newcommand{\R}{\mathbb{R}}
	
	\newcommand{\Z}{\mathbb{Z}}

	\newcommand{\mfX}{\mathfrak{X}}

	\newcommand{\ECS}{manifold with essentially parallel Weyl tensor} 

\begin{document}

	\title{\textsc{On the Full Holonomy of Lorentzian Manifolds with Parallel Weyl Tensor}}
	  
  \author{\textsc{Daniel Schliebner}\thanks{
    The author is funded by the Berlin Mathematical School (BMS).}}
  \date{\Completion}

	\maketitle
	
  \begin{abstract}
  \noindent We compute the full holonomy group of compact Lorentzian manifolds with parallel Weyl tensor, which are neither conformally flat nor locally symmetric, 
  for the case where the fundamental group is contained in a distinguished subgroup $G$ of the isometry group of the universal cover.
  To prove this, we show that every such compact Lorentzian manifold has to be geodesically complete.	  
  Moreover we characterize the identity component of the isometry group for this universal cover and show that $G$ is up to a discrete factor contained in the latter. 
  Concretely, we prove that the identity component of the isometry group is isomorphic to a semidirect product of a subgroup of $\operatorname{SO}(n)$ with 
  the Heisenberg group.
  
  \medskip
  
  
  \end{abstract}

	

	\section{Preliminaries}
\label{section:Pre}
 

We consider Lorentzian manifolds\footnote{In this paper, all manifolds are assumed to be connected and smooth.} $(M^n, g)$ with \textit{essentially parallel Weyl tensor}\footnote{We stress that usually these manifolds are referred to as \textit{essentially conformally symmetric}. Unfortunately this term is a little misleading so we decided to use the term \textit{essentially parallel Weyl tensor} as this describes the properties more carefully.} $W$, i.\,e. which are neither conformally flat ($W = 0$) nor locally symmetric ($\nabla \mathcal R = 0$), where we may assume $n \geq 4$.\footnote{Of course, by \cite[Theorem C]{3} we always have $n \geq 5$ in the Lorentzian case.} This class of manifolds was mainly studied by \textsc{Derdzi\'{n}ski} and \textsc{Roter} in \cite{2,3,5} among others. Admitting a parallel tensor these manifolds need to have special holonomy by the holonomy principle, where we define the \textit{holonomy group} as follows.

\begin{definition}
	\label{Def:Holonomy}
	For any semi-Riemannian manifold $(M^n, g)$ and each $x \in M$ we denote with
	$\Hol_x(M^n,g) := \{ \mathcal P^g_\gamma \in O(T_xM) \ | \ \gamma \in \Omega(x) \} \subset O(T_xM)$	
	the \dStr{holonomy group} of $(M^n, g)$ in $x \in M$ and with
	$\Hol^0_x(M^n,g) := \{ \mathcal P^g_\gamma \in O(T_xM) \ | \ \gamma \in \Omega_0(x) \} \subset O(T_xM)$
	the \dStr{reduced holonomy group} of $(M^n, g)$ in $x \in M$. Here, $\mathcal P^g_\gamma$ denotes the parallel displacement along $\gamma$ belonging to the metric $g$ and
	$\Omega(x)$ denotes the set of piecewise smooth curves closed in $x \in M$ and $\Omega_0(x)$ the subset of curves in $\Omega(x)$ which are null-homotopic.
\end{definition}

Since every pseudo-Riemannian manifold with parallel Weyl tensor admits a totally-isotropic parallel distribution $\mathcal D$, usually referred to as \textit{Olszak distribution}, their holonomy representation leaves $\mathcal D$ invariant. This in turn restricts the holonomy of such a manifold. For instance in the Lorentzian case one especially has $\Hol(M^n, g) \subset (\R^* \times O(n - 2)) \ltimes \R^{n - 2}$. \par	
	Using classification theorems of \textsc{Derdzi\'{n}ski} and \textsc{Roter} in \cite{3} and a result from \cite{1} we prove the following.

\begin{unnumberedtheorem}
	The universal covering $(\widetilde{M}, \widetilde{g})$ of a compact Lorentzian manifold $(M^n, g)$ with essentially parallel 
	Weyl tensor is isometric to $(\R^n, \widetilde{g})$	for a certain complete pp-wave metric $\widetilde{g}$.{\rm \footnote{A Lorentzian metric $h$ on $\R^{1,n-1}$ is 
	called \textit{pp-wave} iff $h(t,s,x) = 2dtds + H(t,x)dt^2 + dx^2$.}} In particular, $(M^n, g)$ must be	itself complete.
\end{unnumberedtheorem}

Using results of \cite{5} and ideas in \cite{7}\footnote{In fact, our computations were motivated by a preliminary version. After finishing this paper, we noticed that in the cited paper, the new Theorem 3 now already implies our Theorem \ref{Thm:Main}.
} this yields the following main theorem concerning the full holonomy group of compact Lorentzian manifolds $(M^n, g)$ with parallel Weyl tensor, whose fundamental group is contained in a certain $(2n - 3)$-dimensional Lie subgroup $G$ of the isometry group $\operatorname{Isom}(\widetilde{M}, \widetilde{g})$ which turns out be isomorphic to $\Z \ltimes \He(n-2)$:

\begin{unnumberedtheorem}
	Let $(M^n, g)$ be a compact Lorentzian manifold with essentially parallel Weyl tensor. Then the reduced holonomy group of $(M^n, g)$ is equal to $\R^{n - 2}$.
	If, moreover, the fundamental group $\pi_1(M)$ is contained in $\Z \ltimes \He(n-2) \subset \operatorname{Isom}(\widetilde{M}, \widetilde{g})$ then
	even $\Hol(M^n, g) = \R^{n - 2}$.
\end{unnumberedtheorem}

We point out that the dimensions $3j + 2$ with $j \geq 1$ are up to now the \textit{only} one in which \textit{compact} pseudo-Riemannian manifolds with parallel Weyl tensor are known to exist, cf. \cite{5}. Away from this restriction, the question that arises is, how restrictive the assumption $\pi_1(M) \subset G$ is. More specifically one could ask, how much larger the full isometry group of $(\widetilde{M}, \widetilde{g})$ can be. Using results of \cite{3} and \cite{8} we will provide a step towards this question by characterizing the identity component of the full isometry group as being isomorphic to $\mathcal S \, {}_{\exp}{\ltimes} \He(n - 2)$ for a certain Lie subgroup $\mathcal S$ of $\operatorname{SO}(n - 2)$. To establish this result we prove in particular that the Lie group $G$ mentioned above is isomorphic to $\Z \ltimes \He(n-2)$. 
%
%
This paper is structured as follows. In Section \ref{section:Lorentz} we present some facts about manifolds with essentially parallel Weyl tensor along with the asserted theorem yielding the characterization in the compact Lorentzian case. Section \ref{section:FH} is then devoted to the computation of the holonomy group and a brief presentation of some results of \cite{5} which are necessary to present the secondly stated theorem. Finally, Section \ref{section:Isom} presents the stated result about the isometry group of the Lorentzian manifolds $(\widetilde{M}, \widetilde{g})$ with essentially parallel Weyl tensor, occurring as the universal cover as above.
	
	\section{On Lorentzian Manifolds with parallel Weyl tensor}
\label{section:Lorentz}


Let $\mathcal D := \bigsqcup_{x \in M} \left\{ u \in T_xM \ | \ g_x(u, \cdot) \wedge W_x(v,v',\cdot,\cdot) = 0 \ \forall v,v' \in T_xM \right\} \subset TM$ denote the Olszak distribution to a pseudo-Riemannian manifold $(M^n, g)$ with parallel Weyl tensor, cf. \cite{2}. It is well known, that the distribution is totally-isotropic and parallel and is either one- or two-dimensional if $M$ is not conformally flat, cf. \cite[Section 2]{2}. For the case that $\dim \mathcal D = 1$, \textsc{Derdzi\'{n}ski} and \textsc{Roter} proved the following.

\begin{theorem}{\bf \cite[Theorem 4.1]{2}}
	\label{Thm:Local}
	Let $(M^n, g)$ be a pseudo-Riemannian manifold with parallel Weyl tensor and $\dim \mathcal D = 1$. 
	Then, every $x \in M^n$ has a connected neighborhood isometric to an open subset of a manifold
	\begin{equation}
		\label{Equation:Local structure of CS with d = 1}
		(I \times \R \times V, \kappa dt^2 + dtds + \Theta).
	\end{equation}
	Thereby, $I \subset \R$ is an open interval and $V$ is a real vector space of dimension $n-2$ with a pseudo-Euklidean inner product 
	$\left\langle \cdot, \cdot\right\rangle$. Furthermore, $t$ and $s$ denote the Cartesian coordinates on the $I \times \R$ factor, $\Theta := \pi^*\left\langle \cdot, \cdot\right\rangle$ for
	$\pi : I \times \R \times V \To V$ and $\kappa : I \times \R \times V \To V$ is defined by $\kappa(t,s,\psi) := f(t)\left\langle \psi,\psi \right\rangle + \left\langle A\psi, \psi \right\rangle$,
	where $f : I \To \R$ is a $C^\infty$-function and $A \in \operatorname{End}(V)$ is a nonzero traceless operator, self-adjoint relative to $\left\langle \cdot, \cdot \right\rangle$. \par
		These manifolds are never conformally flat, and locally symmetric if and only if $f$ is constant.
\end{theorem}

Since the Olszak distribution of a \textit{Lorentzian} manifold with parallel Weyl tensor is always one-dimensional if the manifold is not conformally flat, we may ask, if the manifolds described in the preceding theorem have additional properties. As it turns out, this is indeed the case as the following theorem proves.

\begin{theorem}
	\label{Thm:Complete}
	The simply-connected Lorentzian manifolds $(\R^2 \times V, \widetilde{g})$ with the metric 
	$\widetilde{g} = \kappa dt^2 + dtds + \Theta$ defined in Theorem \ref{Thm:Local} are all geodesically complete.	
\end{theorem}

\begin{proof}
	We use a result of \cite{1} by which the completeness is the case if and only if the maximal solutions $s \longmapsto x(s) \in V$ of the differential equation
	\begin{equation}
		\label{equ:GPW-PDE}
		\frac{\nabla^V\dot x(s)}{ds} = \frac{1}{2}(\grad^V \kappa)(s, x(s))
	\end{equation}
	are defined on the whole $\R$. If we choose an orthonormal basis $v_1, \ldots, v_{n - 2} \in V$ of eigenvectors for the self-adjoint operator 
	$A \in \operatorname{End}(V)$, i.\,e. $Av_i = \lambda_iv_i$, $i = 1, \ldots, n - 2$, then one computes by the definition of the 
	function $\kappa : \R^2 \times V \To \R$ that
	$$
		\grad^V \kappa(s, x(s)) = (2x_1(s)(f(s) + \lambda_1), \ldots, 2x_{n - 2}(s)(f(s) + \lambda_{n - 2}))
	$$
	and hence (\ref{equ:GPW-PDE}) turns into the system of ODEs
	\begin{equation}
		\label{equ:GPW-PDE-2}
		\ddot x_i(s) = (f(s) + \lambda_i)x_i(s).
	\end{equation}
	A solution $s \longmapsto x(s)$ to such a system is always defined for all $s \in \R$, see also the remarks following \cite[Corollary 3.4]{1} and this 
	finishes the proof.
\end{proof}

An obvious consequence of Theorem \ref{Thm:Complete} is that every Lorentzian manifold with essentially parallel Weyl tensor which has a manifold $(\R^2 \times V, \widetilde{g})$ as universal covering, must be itself complete. As it turns out, this is always the case, when $(M^n, g)$ is \textit{compact}\footnote{If one assumes that all leaves of $\mathcal D^\bot$ are complete, this result was already mentioned in \cite{3} by \textsc{Derdzi\'{n}ski} and \textsc{Roter} but without completely proving it. Although this proof is not difficult, the fact that the leaves are always complete motivated us to present the complete proof here.}:

\begin{theorem}
	\label{Thm:UniversalCover}
	Let $(M^n,g)$ be a compact Lorentzian {\ECS}. Then the universal cover
	$\widetilde{M}$, endowed with the induced metric, is isometric to a manifold 
	\begin{equation}
		\label{equ:ECSModel}
		(\R^2 \times V, \kappa dt^2 + dtds + \Theta),
	\end{equation}
  where $V$ is a real vector space of dimension $n-2$ with an Euklidean inner product 
	$\left\langle \cdot, \cdot\right\rangle$. Furthermore, $t$ and $s$ denote the Cartesian coordinates on the $\R^2$ factor, $\Theta := \pi^*\left\langle \cdot, \cdot\right\rangle$ for
	$\pi : \R^2 \times V \To V$ and $\kappa : \R^2 \times V \To V$ is defined by $\kappa(t,s,\psi) := f(t)\left\langle \psi,\psi \right\rangle + \left\langle A\psi, \psi \right\rangle$,
	where $f : \R \To \R$ is a periodic $C^\infty$-function and $A \in \operatorname{End}(V)$ is a nonzero traceless operator, 
	self-adjoint relative to $\left\langle \cdot, \cdot \right\rangle$.
\end{theorem}

Before we prove this assertion, we need to cite some results of \textsc{Derdzi\'{n}ski} in \cite{3}. First we need to introduce a class of manifolds on whose universal cover a one-dimensional factor splits up.

\begin{definition} \label{Def:t-property}
	We say that a manifold $M$ satisfies the \dStr{$t$-property} if its universal covering $\widetilde{M}$ admits a smooth function $t : \widetilde{M} \To \R$ such that
	there exists a manifold $N$ and a diffeomorphism $\psi : \widetilde{M} \To \R \times N$ with the property
	\begin{equation}
		\label{Equation:Universal covering t-property}
		t = \pr_1 \circ \ \psi,
	\end{equation}
	whereby $\pr_1$ is the projection of $\R \times N$ onto $\R$.		
\end{definition}

\noindent With the aid of this definition, \textsc{Derdzi\'{n}ski} proved the following results.\footnote{Of course, although these results are difficult to prove we cite them as lemmata here for a better understanding of the context.}

\begin{lemma}{\bf \cite[Theorem 7.1]{3}}
	\label{lemma:SimplyConnectedECS}
	Let $(M^n,g)$ be a simply connected Lorentzian {\ECS } and let $t : M \To \R$ be a smooth function such that $\grad_g t$ is a global parallel vector field spanning the
	Olszak distribution $\mathcal D$. If, additionally, all leaves of the distribution $\mathcal D^\bot$ are complete and $M$ satisfies the $t$-property, then 
	$(M^n, g)$ is isometric to a manifold occurring in Theorem \ref{Thm:UniversalCover}.
\end{lemma}

\noindent Moreover, one has the following result for the Olszak distribution $\mathcal D$ of $(M^n, g)$.

\begin{lemma}{\bf \cite[Theorem D]{3}}
\label{lemma:Dspanned}
Let $(M^n, g)$ be a Lorentzian {\ECS}. Passing to a two-fold covering manifold, if necessary, we may assume that the Olszak distribution
$\mathcal D$ of $M$ is trivial as a real line bundle, and then $\mathcal D$ is spanned by a global vector field which is even parallel.
\end{lemma}

The final lemma we use is the following which gives a sufficient condition for a manifold $M$ having the $t$-property.

\begin{lemma} 
 \label{Lemma:Universal covering t-property}
	Let $M^n$ be a compact manifold and $\xi \in \Omega^1(M)$ be a nowhere vanishing differential form on $M$. 
	Furthermore, let $\phi : \widetilde{M} \To M$ denote the universal covering map.	
	Then the following holds: If there exists a smooth function $t : \widetilde{M} \To \R$ with $dt = \phi^*\xi$, then there exists a manifold $N$ and a diffeomorphism 
	$\psi : \widetilde{M} \To \R \times N$ with the property $t = \pr_1 \circ \ \psi$, whereby $\pr_1$ is the projection of $\R \times N$ onto $\R$.
\end{lemma}

\begin{proof}
	This is a part of \cite[Lemma 1.2]{3} and can easily be proved.
\end{proof}

\noindent Let us come back to the proof of Theorem \ref{Thm:UniversalCover}. The idea is now, to use the global parallel section $T \in \Gamma(\mathcal D)$ provided by
Lemma \ref{lemma:Dspanned} to show that the leaves of $\mathcal D^\bot$ are always complete with respect to the connection $\nabla^{\mathcal D^\bot}$ induced by the Levi-Civita connection of $g$. Then, applying Lemma \ref{lemma:SimplyConnectedECS} will yield the theorem.

\begin{lemma}
	\label{lemma:CompleteLeaves}
	Let $(M^n, g)$ be a compact Lorentzian manifold with essentially parallel Weyl tensor, whose Olszak distribution $\mathcal D$ is spanned by a global parallel section $T$.
	Then each leaf $L$ of the parallel distribution $\mathcal D^\bot$ is complete w.\,r.\,t. the connection $\nabla^{\mathcal D^\bot}|_L$, 
	induced by the Levi-Civita connection of $g$.
\end{lemma}

\begin{proof}		
	This proof was motivated by \cite{11}. Taking the global section $T \in \Gamma(\mathcal D)$, 
	we may define $Z \in \mfX(M)$ to be the vector field with $g(T,Z) = 1$ and $g(Z,Z) = 0$. Let
	$\mathcal S := \spann\{T,Z\}^{\bot_g} \subset TM$ denote the orthogonal subbundle (usually referred to as the \textit{screen bundle}). We then define by
	$$
		g^R(X,Y) := 
		\begin{cases}
			1, 			& X = Y = T \text{ or } X = Y = Z, \\
			g(X,Y),	& X,Y \in \Gamma(\mathcal S), \\
			0,			& \text{otherwise}.
		\end{cases}
	$$
	a \textit{Riemannian} metric on $M$. We now assert that $\nabla_T^g X \in \Gamma(\mathcal S)$ for all $X \in \Gamma(\mathcal S)$ which 
	indeed is a consequence of Theorem \ref{Thm:Local}. Namely,	let $x \in M$ be arbitrary and $(U, \varphi)$ a chart around $x \in M$ s.\,t. $(U, g|_U)$ 
	is isometric to an open subset of the manifold 
	$(I \times \R \times V, \widetilde{g})$	with $\widetilde{g} = dtds + \kappa dt^2 + \sum_{i = 1}^{n - 2} dx_i^2$,
	occurring in Theorem \ref{Thm:Local} with $t$ resp.\ $s$ denoting the coordinates on $I \subset \R$ resp.\ $\R$ and $x_1, \ldots, x_{n - 2}$ the coordinates on
	$V$. Denoting with $F : \widetilde{U} \subset I \times \R \times V \To U$ the local isometry, cf. \cite[Lemma 5.1]{2}, one can prove
	that $dF_{(t,s,v)}(\partial_s) = \tfrac{1}{2} T(F(t,s,v))$ cf. the proof of \cite[Lemma 5.1]{2}. In particular, any $X \in \Gamma(\mathcal S)$ cannot have
	a $dF_{(t,s,v)}(\partial_t)$-component.
	Now, since	$\nabla_{\partial_s}^{\widetilde{g}} \partial_s = \nabla_{\partial_s}^{\widetilde{g}} \partial_t 	= \nabla_{\partial_s}^{\widetilde{g}} \partial_i = 0$
	for the coordinate vector fields $\partial_s, \partial_t$ and $\partial_i$, $i = 1, \ldots, n - 2$, on $I \times \R \times V$, this proves the assertion. \par
		Let $L$ be a leaf of the foliation induced by the parallel distribution $\mathcal D^\bot$. Taking into account the Koszul formula and 
	$\nabla_T^g X \in \Gamma(\mathcal S)$ for all $X \in \Gamma(\mathcal S)$, 
	one sees that	$\nabla^R|_{L} = \nabla^{\mathcal D^\bot}|_L$, where $\nabla^R|_{L} = \pr_{TL} \circ \nabla^R|_L$ denotes the Levi-Civita connection on $L$ induced by $g^R$.
	The gist now is that the manifold $(L, g^R|_{L \times L})$ occurring as a leaf of the foliated compact Riemannian manifold $(M^n, g^R)$ is always complete, which
	is a well-known fact of foliation theory, see \cite{10} for example. As $\nabla^R|_{L} = \nabla^{\mathcal D^\bot}|_L$ this yields completeness of 
	$\nabla^{\mathcal D^\bot}|_L$.	
\end{proof}

\noindent We can now prove Theorem \ref{Thm:UniversalCover}.

\begin{proof}[Proof of Theorem \ref{Thm:UniversalCover}]	
	By Lemma \ref{lemma:Dspanned}, we may obtain a global parallel section $T \in \Gamma(\mathcal D)$, spanning $\mathcal D$, if $\mathcal D$ is orientable. 
	Otherwise, denoting by 
	${\rm Or}(\mathcal D_x) := \{ \pm \mathcal O_{\mathcal D_x} \}$ the set of the two possible orientations in $x \in M$, we may pass to the two-fold covering manifold
	$\widehat{M} := \bigsqcup_{x \in M} {\rm Or}(\mathcal D_x)$, which is then orientable and connected, as $\mathcal D$ where assumed to be non-orientable. 
	Furthermore, $\widehat{M}$ is compact, since the compact manifold $M$ is finitely covered by $\widehat{M}$.
	Hence, if	$\mathcal D$ is non-orientable, we can pass to the two-fold covering $\widehat{M}$ and as it is connected, it must be covered by the universal covering manifold
	$\widetilde{M}$ as well. Consequently, we may assume w.\,l.\,o.\,g. that $\mathcal D$ is spanned by a global parallel section $T \in \Gamma(\mathcal D)$ (by passing to the
	two-fold covering $\widehat{M}$, if necessary).
	
	\begin{description}
		\item[Step 1: Completeness of the leaves of $\mathcal D^\bot$.] 
		To obtain completeness of the leaves to $\widetilde{\mathcal D}^\bot$ w.\,r.\,t. the connection $\nabla^{\widetilde{\mathcal D}^\bot}$ 
		induced by the Levi-Civita connection of $g$ we take Lemma \ref{lemma:CompleteLeaves} into account. Using the global parallel section as argued before,
		we thus obtain completeness of the leaves of $\mathcal D^\bot$.	
	\end{description}
	
	\begin{description}
		\item[Step 2: $\widetilde{M}$ satisfies the $t$-property.]
	
		Let $\phi : \widetilde{M} \To M$ denote the universal covering map of $M$. Our aim is to apply Lemma \ref{lemma:SimplyConnectedECS} on the
		universal covering $(\widetilde{M}, \widetilde{g})$ of $M$, whereby $\widetilde{g} := \phi^*g$. 		
		Denote by $\mathcal D \subset TM$ and $\widetilde{\mathcal D} \subset T\widetilde{M}$ the the Olszak distributions on $M$ and $\widetilde{M}$, respectively.
		Since the totally-geodesic leaves of $\mathcal D^\bot$ are complete, so are the leaves of $\widetilde{\mathcal D}^\bot$ as $\phi$ is a local isometry.		
		We do now want to apply Lemma \ref{Lemma:Universal covering t-property} in order to obtain that $\widetilde{M}$ satisfies the $t$-property. Thus, the
		global parallel section in Lemma \ref{lemma:SimplyConnectedECS} has to be a pullback of a parallel vector field on $M$, spanning $\mathcal D$. 
		Taking the global parallel section $T \in \Gamma(\mathcal D)$, we define by
		\begin{equation}
			\label{Equation:Proof Universal Cover of Compact Lorentzian ECS no 0}
			\widetilde{T} : x \in \widetilde{M} \longmapsto \widetilde{T}(x) := (\phi^*T)(x) = (d\phi_x)^{-1}(T(\phi(x))) \in \widetilde{\mathcal D}_x
		\end{equation}
		a smooth section, which is parallel, as so is $T$. Indeed, taking an arbitrary path
		$\widetilde{\gamma} : [0,1] \To \widetilde{M}$ with $\widetilde{\gamma}(0) = x_0$, we obtain the formula
		\begin{equation}
			\label{Equation:Proof Universal Cover of Compact Lorentzian ECS no 1}
			\mathcal P_{(\phi \circ \widetilde{\gamma})|_{[0,t)}}^{g} \circ d\phi_{x_0} = 
			d\phi_{\widetilde{\gamma}(t)} \circ \mathcal P_{\widetilde{\gamma}|_{[0,t)}}^{\widetilde{g}},
		\end{equation}
		since $\phi$ is a local isometry. Hence,
		\begin{eqnarray*}
			\mathcal P_{\widetilde{\gamma}|_{[0,t)}}^{\widetilde{g}}(\widetilde{T}(x_0))
				& \stackrel{(\ref{Equation:Proof Universal Cover of Compact Lorentzian ECS no 0})}{=} & 
				\mathcal P_{\widetilde{\gamma}|_{[0,t)}}^{\widetilde{g}} ( (d\phi_{x_0})^{-1}(T(\phi(x_0))) )
				\stackrel{(\ref{Equation:Proof Universal Cover of Compact Lorentzian ECS no 1})}{=}
				(d\phi_{\widetilde{\gamma}(t)})^{-1}( \mathcal P_{(\phi \circ \widetilde{\gamma})|_{[0,t)}}^{g} (T(\phi(x_0))) ) \\
				& \stackrel{(!)}{=} & (d\phi_{\widetilde{\gamma}(t)})^{-1}( T(\phi(\widetilde{\gamma}(t))) )
				\stackrel{(\ref{Equation:Proof Universal Cover of Compact Lorentzian ECS no 0})}{=} \widetilde{T}(\widetilde{\gamma}(t)),
		\end{eqnarray*}
		where $(!)$ holds by parallelity of $T$. Of course, $\widetilde T \in \Gamma(\widetilde{\mathcal D})$ then is a gradient field,
		such that there exists a smooth function $t : \widetilde{M} \To \R$ with $dt = \widetilde{g}(\widetilde T, \cdot) = \phi^*g(T, \cdot)$, where the last equality
		is easily verfied: For $x \in \widetilde{M}$ and $v \in T_x\widetilde{M}$, we obtain
		$$
			dt_x(v) 
			= \widetilde{g}_x(\widetilde{T}(x), v)
			= (\phi^*g)_x(\widetilde{T}(x), v)
			= g_{\phi(x)}(d\phi_x(\widetilde{T}(x)), d\phi_x(v))
			= g_{\phi(x)}(T(\phi(x)), d\phi_x(v)).
		$$
		Thus, by setting $\xi := g(T, \cdot) \in \Omega^1(M)$, we obtain a nowhere vanishing 1-form on $M$, satisfying
		$dt = \phi^*\xi$.	Finally, since $M$ (or, if necessary $\widehat{M}$) is compact, we can apply Lemma \ref{Lemma:Universal covering t-property} and obtain
		a diffeomorphism $\psi : \widetilde{M} \to \R \times N$, whereby $t : \widetilde{M} \To \R$ coincides with the projection $\pr_1 : \R \times N \To \R$, so that
		$\widetilde{M}$ satisfies the $t$-property.
	\end{description}
	
	\noindent Therefore, $\widetilde{M}$, endowed with the pullback metric of $g$, is a simply connected Lorentzian manifold which satisfies the $t$-property and whose
	leaves of the orthogonal complement of the Olszak distribution are all complete. Consequently, applying Lemma \ref{lemma:SimplyConnectedECS} completes
	the proof.
\end{proof}

This yields the following corollary which will enable us to compute the full holonomy group 
of compact Lorentzian manifolds $(M^n, g)$ with essentially parallel Weyl tensor whose fundamental group $\pi_1(M)$ is assumed to be contained 
in a special subgroup of $\operatorname{Isom}(\widetilde{M}, \widetilde{g})$, see Section \ref{section:FH}.

\begin{corollary}
	\label{Cor:CompleteCptLorentz}
	Every compact Lorentzian manifold $(M^n, g)$ with essentially parallel Weyl tensor is isometric to a manifold $\widetilde{M}/\Gamma$, where $\Gamma = \pi_1(M)$ 
	and	$(\widetilde{M}, \widetilde{g})$ is isometric to a manifold occurring in Theorem \ref{Thm:UniversalCover}. In particular,	$(M^n, g)$ is geodesically complete due to
	Theorem \ref{Thm:Complete}.
\end{corollary}
	
	\section{On the Full Holonomy Group}
\label{section:FH}

The aim of this section is to compute the full holonomy group of certain compact Lorentzian manifolds $(M^n, g)$ with essentially parallel Weyl tensor. Due to Theorem \ref{Thm:UniversalCover} such a manifold is isometric to a quotient $\widetilde{M}/\Gamma$ where $(\widetilde M, \widetilde{g}) = (\R^2 \times V, dtds + \kappa dt^2 + \Theta)$ is the Lorentzian universal covering manifold of $M$. If we consider the Lorentzian universal covering map
$$
	(\R^2 \times V, \widetilde{g}) \To (\widetilde{M}/\Gamma, [\widetilde{g}]) \simeq (M^n, g)
$$
the idea is now to use the following lemma in order to compute the full holonomy group of the manifold $\widetilde{M}/\Gamma \simeq M$ where we require $\pi_1(M)$ to be contained in a special subgroup $G \subset \operatorname{Isom}(\widetilde{M}, \widetilde{g})$ (see below).

\begin{lemma}
	\label{lemma:HolOfCovering}
	Let $\pi : (\widetilde M, \widetilde{g}) \To (M^n, g)$ be a semi-Riemannian covering map. If $\widetilde{\gamma}$ is a lift of a loop $\gamma : [0,1] \To M$ and 
	$\mathcal P^{\widetilde{g}}_{\widetilde{\gamma}}$ resp. $\mathcal P^{g}_{\gamma}$ denote the parallel displacements along $\widetilde{\gamma}$ resp. $\gamma$, then
	it holds
	$$
		\mathcal P^{g}_{\gamma}(v) = d\pi_{\widetilde{\gamma}(1)}(\mathcal P^{\widetilde{g}}_{\widetilde{\gamma}}(\widetilde{v})),
	$$
	where $\widetilde{v} \in T_{\widetilde{\gamma}(0)}\widetilde M$ and $d\pi_{\widetilde{\gamma}(0)}(\widetilde{v}) = v$.
\end{lemma}

\noindent Using this together with the fact that the map
\begin{equation}
	\label{Equ:SurjHol}
	\pi_1(M, p) \ni [\gamma] \longmapsto [\mathcal P^g_{\gamma}] \in \Hol_p(M, g)/\Hol_p^0(M, g) \ \text{ with } p \in M
\end{equation}
is a surjection (cf. \cite[10.15]{4}), we can compute the full holonomy group $\Hol_p(M, g)$. However, for concrete calculations we need to know how the subgroups $\pi_1(M) = \Gamma \subset \operatorname{Isom}(\widetilde{M}, \widetilde{g})$, producing a compact Lorentzian quotient manifold with essentially parallel Weyl tensor, look like. This is, were the restriction onto the fundamental group of $M$ comes into play. Namely, \textsc{Derdzi\'{n}ski} investigated in \cite{5} discrete subgroups 
$\Gamma \subset G \subset \operatorname{Isom}(\widetilde{M}, \widetilde{g})$ for a certain group of isometries $G$ that produce compact quotients. Note that these examples
are also complete by Theorem \ref{Thm:Complete} and thus fit into our setting. We use their explicit characterization (cf. \cite[Theorem 6.1]{5}) of such groups $\Gamma$ in order to compute the full holonomy group $\Hol(M, g)$ for those manifolds with $\pi_1(M) = \Gamma \subset G$. This is done by the above-mentioned observations together with noting that we are already aware of the reduced holonomy group $\Hol^0(M, g)$ since this is isomorphic to the full holonomy group $\Hol(\widetilde{M}, \widetilde{g})$ of the universal cover, cf. \cite[Theorem 5.4]{6}. In \cite[Example 5.5]{6} there is also shown that for $(\widetilde{M}^n, \widetilde{g}) = (\R^2 \times V, dtds + \kappa dt^2 + \Theta)$ we have
$$
	\Hol(\widetilde{M}^n, \widetilde{g}) = \R^{n - 2}
$$
since $(V, \left\langle \cdot, \cdot \right\rangle)$ is flat. Thus we have just proven the following.

\begin{lemma}
	\label{lemma:RestrHol}
	Let $(M^n, g)$ be a compact Lorentzian manifold with essentially parallel Weyl tensor. Then the reduced holonomy group is 
	given by $\Hol^0(M^n, g) = \R^{n-2}$.
\end{lemma}

Before we proceed to compute the holonomy group $\Hol_p(M, g)$ let us briefly explain the construction of \textsc{Derdzi\'{n}ski}. For a complete description we refer the reader to the original paper \cite{5}. As already mentioned, \textsc{Derdzi\'{n}ski} considered suitable discrete subgroups $\Gamma \subset G \subset \operatorname{Isom}(\widetilde{M}, \widetilde{g})$ which produce \textit{compact} quotient manifolds $\widetilde{M}/\Gamma$ with essentially parallel Weyl tensor. Namely, the subgroup $G \subset \operatorname{Isom}(\widetilde{M}, \widetilde{g})$ is therein defined as follows. \par
	Let $f, A, V$ and $\left\langle \cdot, \cdot \right\rangle$ denote the objects occurring in the classification theorem (Theorem \ref{Thm:UniversalCover}).\footnote{Note that in \cite{5}, $f \in C^\infty(\R)$ is assumed to be periodic of period $p$. This subtlety is however no restriction here as this is for the Lorentzian case due to \cite[Remark 5.1]{3}.} Given the solution space 
\begin{equation}
	\label{Equ:E}
	\mathcal E := \{ u : \R \To V \text{ smooth} \ | \ \ddot u(t) = f(t)u(t) + Au(t) \},
\end{equation}
the non-degenerate skew-symmetric bilinear form
$$
	\Omega(u_1, u_2) := \left\langle \dot u_1, u_2 \right\rangle - \left\langle u_1, \dot u_2 \right\rangle
$$
which is constant for all $t \in \R$ since $A \in \operatorname{End}(V)$ is self-adjoint relative to $\left\langle \cdot, \cdot \right\rangle$, and the linear isomorphism $T : \mathcal E \To \mathcal E$ with $(Tu)(t) := u(t - p)$ (where $p \in \R$ denotes the period of $f \in C^\infty(\R)$), one defines
\begin{equation}
	\label{Equ:G}
	G := \Z \times \R \times \mathcal E
\end{equation}
with $g_1 \cdot g_2$ for $g_i := (k_i, x_i, u_i)$ defined through the formula
\begin{equation}
	\label{Equ:Gop}
	g_1 \cdot g_2 := (k_1 + k_2, x_1 + x_2 - \Omega(u_1, T^{k_1}u_2), T^{-k_2}u_1 + u_2).
\end{equation}
For $g = (k,x,u)$ and $m = (t,s,v) \in \widetilde{M} = \R \times \R \times V$, the action of $G$ on $\widetilde{M}$ then is given through
\begin{equation}
	\label{Equ:Gact}
	g \cdot m := (t + kp, s + x - \left\langle \dot u(t), 2v + u(t) \right\rangle, v + u(t)).
\end{equation}
One easily verifies that each $F_g : m \in \widetilde{M} \To g \cdot m \in \widetilde{M}$ is an isometry for the metrics defined in Theorem \ref{Thm:UniversalCover}, i.\,e. $G \subset \operatorname{Isom}(\widetilde{M}, \widetilde{g})$ as desired. \par
	As the non-degenerate skew-symmetric bilinear form $\Omega \in \Lambda^2\mathcal E^*$ is constant, the pair 
	$(\mathcal E, \Omega)$ defines a symplectic vector space of dimension 
$2(n - 2)$. To every such vector space we can associate a Heisenberg group $\He(\mathcal E, \Omega)$, cf. \cite[Section I.3]{9}, if we endow $\He(\mathcal E, \Omega) := \R \times \mathcal E$ with the group structure $g_1g_2 := (t_1 + t_2 + \Omega(u_1,u_2), u_1 + u_2)$ for $g_i = (t_i, u_i)$. By fixing a Darboux basis in $\mathcal E$, this group is isomorphic the canonical Heisenberg group $\He(n - 2)$ which in matrix representation is defined as
$$
	\He(n - 2) := 
	\left\{\left.\begin{pmatrix}
	1 & a^T & c \\
	0 & \mathbb{I} & b \\
	0 & 0 & 1
	\end{pmatrix} \ \right| \ a,b \in \R^{n - 2}, \ c \in \R \right\}.
$$
\noindent Moreover, we define for $B : \R \To \operatorname{End}(V)$ with $\dot B + B^2 = f \cdot \mathbb{I} + A$ by
\begin{equation}
	\label{Equ:L}
	\mathcal L := \{ u : \R \To V \text{ smooth} \ | \ \dot u(t) = B(t)u(t) \} \subset \mathcal E
\end{equation}
an $(n - 2)$-dimensional subspace of $\mathcal E$. In the classification theorem \cite[Theorem 6.1]{5} it now turns out that each subgroup $\Gamma \subset G$ producing a compact quotient manifold $\widetilde{M}/\Gamma$ implies the existence of a normal subgroup $\Sigma \subset \R \times \mathcal L$ which is a lattice and defined by
$\Sigma := \Gamma \cap \ker \Pi$ with $\Pi : G \twoheadrightarrow \Z$ being the surjective homomorphism defined by $\Pi(k,q,u) := k$, see \cite[Section 4]{5} and \cite[Theorem 6.1]{5} for details. \par
	Using this, we can describe $G$ and $\Gamma$ more explicitly in terms of the groups $\Z$, $\Sigma$ and $\He(n-2)$. Namely, since $\Z$ is a free group and 
	$\phi : (t,u) \in \R \times \mathcal E \cong \{0\} \times \R \times \mathcal E \subset G \longmapsto (-t, u) \in \He(\mathcal E, \Omega)$ 
	is a Lie group isomorphism, the short exact sequences
\begin{eqnarray}	
0 \To \R \times \mathcal E \cong_\phi \He(\mathcal E, \Omega) \cong \He(n-2) \xhookrightarrow{ \ \ \iota \ \ } G \xtwoheadrightarrow{\Pi} \Z \To 0 \label{equ:SplitExact2} \\
0 \To \Sigma \xhookrightarrow{ \ \ \iota \ \ } \Gamma \xtwoheadrightarrow{\Pi} \Z \To 0 \label{equ:SplitExact}
\end{eqnarray}
split. Thus, $G \cong \Z \ltimes \He(n - 2)$ resp. $\Gamma \cong \Z \ltimes \Sigma$ by which we obtain the following lemma.

\begin{lemma}
	\label{lemma:Splits}
	For $G$ and $\Sigma$ as above it holds:
	\begin{itemize}
		\item[(i)] $G \cong \Z \ltimes \He(n - 2)$.
		\item[(ii)] $\Gamma \cong \Z \ltimes \Sigma \cong \Z \cdot \Z^{n - 1}$.
	\end{itemize}	
\end{lemma}

We are now in the position to state the main theorem which describes the full holonomy group of compact Lorentzian manifolds with essentially parallel Weyl tensor whose fundamental group is contained in $G$.

\begin{theorem}
	\label{Thm:Main}
	Let $(M^n, g)$ be a compact Lorentzian manifold with essentially parallel Weyl tensor. Then the reduced holonomy group of $(M^n, g)$ is equal to $\R^{n - 2}$.
	If, moreover, the fundamental group $\pi_1(M)$ is contained in the subgroup $G \subset \operatorname{Isom}(\widetilde{M}, \widetilde{g})$ of the isometry group 
	of the universal cover $(\widetilde{M}, \widetilde{g})$, cf. (\ref{Equ:G}), then it holds
	$$
		\Hol(M^n, g) = \R^{n - 2}
	$$
	for the full holonomy group of $(M^n, g)$.
\end{theorem}

\begin{proof}
	Let us denote the epimorphism (\ref{Equ:SurjHol}) by $\Phi$ and consider the short exact sequence
	\begin{equation}
		\label{equ:ProofMain1}
		0 \To \Hol^0(M,g) \xhookrightarrow{ \ \ i \ \ } \Hol(M,g) \xtwoheadrightarrow{ \ p \ } \Hol(M,g)/\Hol^0(M,g) \stackrel{\Phi^*}{\cong} \pi_1(M)/\ker \Phi \To 0
	\end{equation}
	where $\Phi^* : \pi_1(M)/\ker \Phi \To \Hol(M,g)/\Hol^0(M,g)$ is the isomorphism induced by the group homomorphism theorem. 
	If we assume that the quotient group $H := \pi_1(M)/\ker \Phi$ is free, the sequence (\ref{equ:ProofMain1}) splits, whence we obtain
	$\Hol(M,g) \cong H \ltimes \Hol^0(M,g)$. The idea now is to calculate the parallel transport $\mathcal P_{\gamma}^g$ along the generators of
	$[\gamma] \in \pi_1(M) = \Gamma \cong \Z \ltimes \Sigma$ using Lemma \ref{lemma:HolOfCovering} and to study, 
	in which cases they produce non-trivial elements in $\Hol(M,g)/\Hol^0(M,g)$. \par
		We proceed to construct such generators for each $\sigma \in \Gamma$. For every $\sigma \in \Gamma \cong \Z \ltimes \Sigma$ with		
	$\sigma = (k, \sigma_\Sigma) = (k, (r, w))$, $k \in \Z$, $(r,w) \in \Sigma$, we define by 
	\begin{eqnarray}
		\label{equ:ProofMain2a}
		\widetilde{\gamma}_{\sigma_\Sigma}(t) & := & (0, t \cdot (r - \bigl\langle \dot w(0), w(0) \bigr\rangle), t \cdot w(0)), \\
		\label{equ:ProofMain2b}
		\widetilde{\gamma}_{k}(t) & := & (t \cdot kp, 0, 0),
	\end{eqnarray}
	curves $\widetilde{\gamma}_{\sigma_\Sigma}, \widetilde{\gamma}_k : [0,1] \To \widetilde M$. Indeed, these curves induce generators
	$\gamma_{\sigma_\Sigma} := \pi \circ \widetilde{\gamma}_{\sigma_\Sigma}$ and $\gamma_k := \pi \circ \widetilde{\gamma}_k$ of $\pi_1(M)$, where	
	$\pi : (\widetilde M, \widetilde g) \To (\widetilde M/\Gamma, g)$ denotes the universal covering map. Namely each such curve in closed since
	$\gamma_{\sigma_\Sigma}(0) = \gamma_k(0) = [(0,0,0)]$ and as $\sigma_{\Sigma}^{-1} = (0, -r, -w)$ it holds
	$$
		(0, -r, -w) \cdot \widetilde{\gamma}_{\sigma_\Sigma}(1) = (0,0,0) = (-k, 0, 0) \cdot \widetilde{\gamma}_k(1)
	$$
	by the definition of the action of $G$ on $\widetilde{M}$, cf.	(\ref{Equ:Gact}). Moreover, each two such curves cannot be homotopic if $\sigma_1 \neq \sigma_2$ 
	which is immediate by the proper discontinuous action of $\Gamma$ on $\widetilde{M}$.	\par
		Following the notations in Theorem \ref{Thm:UniversalCover} (with 
	$\partial_i := \frac{\partial}{\partial x_i}$ s.\,t. $\left\langle \partial_i, \partial_j \right\rangle = \delta_{ij}$) we obtain for the Levi-Civita 
	connection $\widetilde{\nabla}$ of $\widetilde{g}$ that
	$$
		\widetilde{\nabla}_{\partial_i} \partial_t = \widetilde{\nabla}_{\partial_t} \partial_i = \partial_i(\kappa)\partial_s
		\text{ and }
		\widetilde{\nabla}_{\partial_t} \partial_t = \partial_t(\kappa)\partial_s - \tfrac{1}{2}\partial^i(\kappa)\partial_i
	$$
	are the only non-vanishing terms. Therefore, $\partial_s$, $\partial_i$ and $\partial_t - \kappa \partial_s$ are parallel along 
	every $\widetilde{\gamma}_{\sigma_\Sigma}$. In particular, they are also parallel along every $\widetilde{\gamma}_k$ as the derivatives $\partial_i(\kappa)$ and
	$\partial_t(\kappa)$ vanish when evaluated in the points $\widetilde{\gamma}_k(t)$ for all $t \in \R$. \par
		Hence, we can compute the parallel transports along the curves $\gamma_{\sigma_\Sigma}$ and $\gamma_k$ using Lemma \ref{lemma:HolOfCovering} easily. Let us fix the basis
	\begin{equation}
		\label{equ:ProofMain3}
		S   := d\pi_{(0,0,0)}(\partial_s), \
		E_i := d\pi_{(0,0,0)}(\partial_i), \
		T := d\pi_{(0,0,0)}(2\{\partial_t - \kappa(0)\partial_s\})
	\end{equation}
	in $T_{\pi(0,0,0)}M$. As $\partial_s$ and $\partial_i$ are parallel along $\gamma_{\sigma_\Sigma}$ we obtain
	\begin{eqnarray*}
		\mathcal P_{\sigma_\Sigma}^g(S) & = & d\pi_{\widetilde{\gamma}_\Sigma(1)}(\partial_s)
		= \left.\frac{d}{da}\pi(0,a + r - \bigl\langle \dot w(0), w(0) \bigr\rangle, w(0))\right|_{a = 0} \\
		& = & d\pi_{(0,0,0)}(\partial_s) = S, \\
		\mathcal P_{\sigma_\Sigma}^g(E_i) & = & d\pi_{(0,0,0)}(\partial_i) = E_i + 2\left\langle \dot w(0), \partial_i \right\rangle S.
	\end{eqnarray*}
	As $\partial_t - \kappa \partial_s$ is parallel along $\gamma_{\sigma_\Sigma}$ too, we obtain after an easy computation that
	\begin{eqnarray*}
		\mathcal P_{\sigma_\Sigma}^g(T) & = & 2 d\pi_{\widetilde{\gamma}_\Sigma(1)}(\partial_t - \kappa(\widetilde{\gamma}_\Sigma(1))\partial_s)
		 = T - \Xi_1(k,w) S - \sum_{i = 1}^{n - 2} \Xi^i_2(w_i)E_i,
	\end{eqnarray*}
	where $\Xi_1$ and $\Xi_2$ denote the maps
	\begin{eqnarray*}
		\Xi_1(w) & := &   2\bigl\langle \dot w(0), \dot w(0) \bigr\rangle, \\
		\Xi^i_2(w) & := & 2\bigl\langle \dot w(0), \partial_i \bigr\rangle.
	\end{eqnarray*}
	Therefore, on the one hand, the only non-trivial parallel transports along $\gamma_{\sigma_\Sigma}$ with $\sigma_\Sigma = (r,w)$ are, 
	written in the basis (\ref{equ:ProofMain3}),
	$$
		\mathcal P^g_{\gamma_{\sigma_\Sigma}} = 
		\begin{pmatrix}
			1 & \Xi^1_2(w) & \ldots & \Xi^{n - 2}_2(w) & -\Xi_1(w) \ \ \ \ \\
			0 & 1 & \ldots & 0 & -\Xi^1_2(w) \ \ \ \ \\
			\vdots	& \vdots & \ddots &   & \vdots \\
			0 & 0 & \ldots & 1 & -\Xi^{n - 2}_2(w) \\
			0 & 0 & \ldots & 0 & 1
		\end{pmatrix}
	$$
	with $w \neq 0$ but on the other hand we have that $\mathcal P^g_{\gamma_{\sigma_\Sigma}} \in \Hol^0(M^n, g) = \R^n$ for all $w \in \mathcal L$. \par
		If we consider the parallel transports along the curves $\gamma_k$, we obtain by the same computations that
	$\mathcal P^g_{\gamma_{k}}(S) = S$, $\mathcal P^g_{\gamma_{k}}(E_i) = E_i$ and $\mathcal P^g_{\gamma_{k}}(T) = T$, i.\,e. that 
	$\mathcal P^g_{\gamma_{k}} = {\rm id}_{T_{\pi(0,0,0)}M}$.
	%
	%
	%
	%
	We thus infer $\ker\Phi = \Gamma$ and if we set $p = \pi(0,0,0)$ this yields
	$$
		\Hol_p(M^n,g)/\Hol_p^0(M^n,g) \cong \pi_1(M, p)/\ker \Phi = \{e\}.
	$$
	With respect to the argumentations at the beginning this proves the theorem.
\end{proof}
	
	\section{On the Isometry Group}
\label{section:Isom}

The aim of this section is to study how ``strong" the restriction on the fundamental group made in Theorem \ref{Thm:Main} is. Indeed, if we still denote with
$(\widetilde{M}^n, \widetilde{g})$ the manifold $(\R^2 \times V, dtds + \kappa dt^2 + \Theta)$ as in Theorem \ref{Thm:UniversalCover}, the fact that the group $\R \times \mathcal E$ as subgroup of $\operatorname{Isom}(\widetilde{M}^n, \widetilde{g})$ is isomorphic to the $(n-2)$-dimensional Heisenberg group $\He(n-2)$, cf. Lemma \ref{lemma:Splits} in Section \ref{section:FH}, is not surprising. Namely, \textsc{Blau} et al.\ \cite{8} proved that for pp-wave metrics with $H(t,x) = \sum_{i = 1}^{n - 2}K(t)x_i^2$, the Lie algebra $\LAkill(\widetilde{M}^n, \widetilde{g})$ is equal to the Lie algebra $\he(n-2)$ of $\He(n-2)$ if $K$ does not fulfill special properties. As $\He(n-2)$ is simply-connected we infer for this case that $\Isom^0(\widetilde{M}^n, \widetilde{g}) \cong \He(n-2)$. Therefore the restriction on the fundamental group $\pi_1(M) = \Gamma$ made in Theorem \ref{Thm:Main} turns out to be equivalent to requiring that $\pi_1(M)$ is contained in $\Z \ltimes \Isom^0(\widetilde{M}^n, \widetilde{g})$. However, if the smooth function $\kappa \in \R^2 \times V \To \R$ has additional properties, the identity component $\Isom^0(\widetilde{M}^n, \widetilde{g})$ might be larger. \par
	The present section should now give a complete answer to the question, how the identity	component $\Isom^0(\widetilde{M}^n, \widetilde{g})$ could look like. In particular we will see that the dimension $d$ (as manifold) of $\Isom^0(\widetilde{M}^n, \widetilde{g})$ (and consequently of $\Isom(\widetilde{M}^n, \widetilde{g})$) can be pretty large, i.\,e. we will prove that for particular $A \in \operatorname{End}(V)$ we have $d = \frac{n(n+1)}{2} - (2n - 3)$. We begin with the asserted description of $\Isom^0(\widetilde{M}^n, \widetilde{g})$.
	
\begin{theorem}
	\label{Thm:Isom0}
	Let $(\widetilde{M}^n, \widetilde{g})$ denote a Lorentzian manifold with essentially parallel Weyl tensor as in Theorem \ref{Thm:UniversalCover}, 
	cf. (\ref{equ:ECSModel}). Then the identity
	component $\Isom^0(\widetilde{M}^n, \widetilde{g})$ of the isometry group of $(\widetilde{M}^n, \widetilde{g})$ is isomorphic to 
	$\mathcal S \, {}_{\exp}{\ltimes} \He(n - 2)$, where $\mathcal S \subset \operatorname{SO}(n - 2)$ is a connected Lie subgroup of $\operatorname{SO}(n - 2)$ 
	with Lie algebra $\mathfrak{s} := \spann\{ F \in \mathfrak{so}(n-2) \ | \ [A,F] = 0 \}${\rm \footnote{Note that $\mathfrak s$ is \textit{anti}-isomorphic to a
	subalgebra of $\mathfrak{so}(n-2)$, i.\,e. $[X,Y]_\mathfrak s = -[X,Y]_{\mathfrak{so}(n-2)}$.}}
	which is non-trivial if and only if $A \in \operatorname{End}(V)$ has at least one eigenspace of dimension greater than one.
\end{theorem}

\begin{proof}	
	In \cite[Section 2.3]{8} it is shown that there always exist $2n - 3$ distinct Killing vector fields $E_1, \ldots, E_{n - 2}, E_1^*, \ldots, E_{n - 2}^*, Z$ for $\widetilde{g}$ which span $\he(n-2)$. Note that in our notation the matrix $A(x^+)$ in \cite{8} equals $(f(t) + \lambda_i)\delta_{ij}$ where $\lambda_i$, $i = 1, \ldots, n - 2$, denote the eigenvalues of $A \in \operatorname{End}(V)$. Indeed, choosing an orthonormal frame in $V$ consisting of eigenvectors $(v_1, \ldots, v_{n - 2})$ and $x_i(t,s,v) := r_i \in \R$, $v = \sum_{i} r_iv_i$, our $\widetilde{g}$ becomes
$$
	\widetilde{g} = dtds + \sum_{i = 1}^{n - 2} (f(t) + \lambda_i)x_i^2 dt^2 + \sum_{i = 1}^{n - 2}dx_i^2.
$$
Concretely, these Killing vector fields $E_1, \ldots, E_{n - 2}, E_1^*, \ldots, E_{n - 2}^*, Z$ are given through
\begin{eqnarray}
	E_i   & = & \sum\nolimits_{k = 1}^{n - 2} (\xi_{i,k} \partial_k - \dot\xi_{i,k}x_k \partial_s), \label{Equ:ProofIsom03a} \\
	E_i^* & = & \sum\nolimits_{k = 1}^{n - 2} (\xi_{i,k}^* \partial_k - \dot\xi_{i,k}^* x_k \partial_s), \label{Equ:ProofIsom03b} \\
	Z     & = & \partial_s, \label{Equ:ProofIsom03c}
\end{eqnarray}
where $(\xi_1, \ldots, \xi_{n - 2},\xi_1^*, \ldots, \xi_{n - 2}^*)$ is the basis of $\mathcal E$ with $\xi_{i,k}(0) = \delta_{ik}$, $\dot\xi_{i,k}(0) = 0$,
$\xi^*_{i,k}(0) = 0$ and $\dot\xi^*_{i,k}(0) = \delta_{ik}$, cf. \cite[Section 2.3]{8}. \par	
	However, there may exist additional Killing vector fields. Namely, this occurs in three distinct cases. The first case produces one additional Killing vector, iff the matrix $K(t) := (f(t) + \lambda_i)\delta_{ij}$ is degenerate for all $t \in \R$ which in our case obviously cannot occur as $f \in C^\infty(\R)$ is non-constant. Next, the authors of \cite{8} determine those cases in which they get additional Killing vector fields with $\partial_t$-component (which they refer to as \textit{homogeneous} pp-waves). As it turns out, these require $K(t)$ to be of the form 
$D(t) = \exp(t F)D_0\exp(-tF)$ or $D(t) = \frac{1}{t^2} \cdot \exp(\ln(t) F)D_0\exp(-\ln(t)F)$ for skew-symmetric $F$ and symmetric $D_0$. But in our case, $D(t)$ cannot be of this form since $\trace K(t) = \trace (f(t) + \lambda_i)\delta_{ij} = (n - 2)f(t)$ but $\trace D(t) = \trace D_0 \equiv const$ or
$\trace D(t) = \frac{1}{t^2} \cdot const$ while $f$ must be periodic. \par
	The remaining case in which additional Killing vector fields can occur is the following. Considering for $F \in \mathfrak{so}(n - 2)$ the equation
\begin{equation}
	\label{Equ:ProofIsom01}
	K(t) \cdot F - F \cdot K(t) \equiv 0,
\end{equation}
it is proven in \cite{8} that each such $F$ generates an additional Killing vector field
\begin{equation}
	\label{Equ:ProofIsom02}
	X_F := \sum_{i,j = 1}^{n - 2} F_{ji}x^i\partial_j.
\end{equation}
Writing down (\ref{Equ:ProofIsom01}) componentwise one sees that $F_{ij} \neq 0$ requires $\lambda_i = \lambda_j$ for the eigenvalues $\lambda_i,\lambda_j \in \R$ of $A$. Thus, non-trivial solutions of (\ref{Equ:ProofIsom01}) occur if and only if $A \in \operatorname{End}(V)$ has at least one eigenspace of dimension greater than one. Taking into account the generated Killing fields (\ref{Equ:ProofIsom02}) we obtain for each such $F \in \mathfrak{so}(n - 2)$ solving (\ref{Equ:ProofIsom01}) a one-parameter subgroup of $\operatorname{SO}(n - 2)$. Namely, the integral curves $\gamma_x^F : \R \To V$ through $x \in V$ of such an $X_F$ satisfy the differential equation
$\dot \gamma_x^F(t) = F \cdot \gamma_x^F(t)$ with $\gamma(0) = x$ whose solution if given through $\gamma(t) = \exp(t \cdot F) \cdot x$ defined on the whole real line. Thus every Killing vector field $X_F$ gives rise to a one-parameter subgroup $\{\Phi^F_\tau\}_{\tau \in \R}$ of $\Isom(\widetilde{M}^n, \widetilde{g})$ via
$$
	\Phi_\tau : (t,s,v) \in \widetilde{M} = \R^2 \times V \longmapsto (t,s,\exp(\tau F)v) \in \widetilde{M}.
$$
Indeed, every $\exp(\tau F) \in \operatorname{SO}(n - 2)$ commutes with $A$ as can be seen by differentiating the equation (which holds since $AF = FA$)
$$
	\exp(sA)\exp(\tau F) = \exp(sA + \tau F) = \exp(\tau F)\exp(sA)
$$
in $s = 0$. Therefore, $\kappa \circ \Phi_\tau = \kappa$ and each $\Phi_\tau$ is an isometry for $\widetilde{g}$. The one-parameter groups now generate the Lie group 
$\mathcal S$ stated in the theorem. Taking account into the commutator relations of elements in $\he(n-2)$ with elements in $\mathfrak{s}$, namely
$$
	[X_F, E_i] = \sum_{\ell = 1}^{n - 2} F_{i\ell}E_\ell, \ \ \ 
	[X_F, E_i^*] = \sum_{\ell = 1}^{n - 2} F_{i\ell}E_\ell^* \ \text{ and } \
	[X_F, Z]     = 0,
$$
cf. (\ref{Equ:ProofIsom03a}) -- (\ref{Equ:ProofIsom03c}) and (\ref{Equ:ProofIsom02}), and considering the homomorphism 
$\pi : \mathcal S \To \operatorname{Aut}(\He(n-2))$ defined through 
$$
	\pi(\exp(F))
	\begin{pmatrix}
		1 & a^T & c \\
		0 & \mathbb{I}   & b \\
		0 & 0 & 1
	\end{pmatrix}
	 :=
	\begin{pmatrix}
		1 & (e^Fa)^T & c \\
		0 & \mathbb{I}   & e^Fb \\
		0 & 0 & 1
	\end{pmatrix}
$$
this yields $\Isom^0(\widetilde{M}, \widetilde{g}) \cong \mathcal S \, {}_{\pi}{\ltimes} \He(n - 2)$. Finally, note that the first set of Killing vector fields
treated as elements of $\mathfrak{he}(n-2)$ are precisely those vector fields with coefficients $\xi_{i,k} \in C^\infty(\R)$ s.\,t. 
$\xi_i = (\xi_{i,1}, \ldots, \xi_{i,n - 2}) \in \mathcal E$. As solutions $\xi \in \mathcal E$ and the additional flows $\Phi_\tau$ as above are all defined on the whole real line this in particular shows that $\mathfrak{kill}(\widetilde{M}, \widetilde{g}) \cong \mathfrak{kill}_c(\widetilde{M}, \widetilde{g})$ with $\mathfrak{kill}_c(\widetilde{M}, \widetilde{g})$ denoting the Lie subalgebra consisting of complete Killing vector fields. As $\mathfrak{kill}_c$ is the Lie algebra corresponding to $\Isom^0(\widetilde{M}, \widetilde{g})$ this completes the proof.
\end{proof}

Let us give an explicit example with $n = 5$ for which we can explicitly compute the Lie group $\mathcal S$ occurring in Theorem \ref{Thm:Isom0}. As it turns out, this example is a special case in $n = 5$ for the subsequent corollary on the maximal dimension of $\Isom^0(\widetilde{M}^n, \widetilde{g})$.

\begin{example}
Let $A \in \End(V)$, $\dim V = 3$, be given through 
$$
	A = \begin{pmatrix}\lambda & 0 & 0 \\ 0 & \lambda & 0 \\ 0 & 0 & \mu \end{pmatrix}
$$
for $\lambda,\mu \in \R \setminus\{0\}$ with $2\lambda + \mu = 0$. Then, up to multiplication by scalars, 
$$
	F = \begin{pmatrix}0 & 1 & 0 \\ -1 & 0 & 0 \\ 0 & 0 & 0 \end{pmatrix}
$$
is the only non-trivial possible solution to (\ref{Equ:ProofIsom01}). Using \textit{Rodrigues's formula} for the matrix exponential
$\exp : \mathfrak{so}(3) \To \operatorname{SO}(3)$, given through
$$
	\exp(X) 
	:= \exp\begin{pmatrix}0 & -c & b \\ c & 0 & -a \\ -b & a & 0 \end{pmatrix} 
	 = \mathbb{I} + \frac{\sin \theta}{\theta}X + \frac{(1 - \cos \theta)}{\theta^2}X^2
$$
where $\theta := \sqrt{a^2 + b^2 + c^2}$, we conclude $\ker(t \mapsto \exp(tF)) = 2\pi\Z$. Hence 
$\{\Phi_t\} \cong \R/2\pi\Z = S^1$ and thus $\Isom^0(\widetilde{M}^n, \widetilde{g}) = S^1 \ltimes \He(3)$.
\end{example}

Indeed, if we generalize the preceding example, we may state the following corollary from Theorem \ref{Thm:Isom0}.

\begin{corollary}
	\label{cor:ev}
	Let $\dim V = n - 2$ and $A \in \End(V)$ be a traceless self-adjoint operator such that $V = E_1 \oplus E_2$ is a splitting of $V$, where
	$E_1$, $E_2$ are the eigenspaces of $A$ such that $\dim E_1 = n - 3$, $\dim E_2 = 1$. Moreover, let $(\widetilde{M}^n, \widetilde{g})$ be constructed as in 
	Theorem \ref{Thm:UniversalCover} using this operator $A$.	Then $\Isom^0(\widetilde{M}^n, \widetilde{g})$ is a Lie subgroup of the full isometry group of 
	dimension $\frac{n(n+1)}{2} - (2n - 3)$.
\end{corollary}

\begin{proof}
	The only restriction in the choice of an $F \in \mathfrak{so}(n - 2)$ satisfying (\ref{Equ:ProofIsom01}) is given through the distinct eigenvalues for the two
	eigenspaces of $A$, i.\,e. $F : V = E_1 \oplus E_2 \To E_1 \oplus E_2$ has to leave the eigenspaces invariant. As $F$ is skew-adjoint and $\dim E_1 = 1$
	we obtain for the Lie algebra $\mathfrak{s} := \spann\{ F \in \mathfrak{so}(n-2) \ | \ [A,F] = 0 \}$ that
	$\dim\mathfrak{s} = \dim \mathfrak{so}(n - 3) = \frac{(n-3)(n-4)}{2}$. Hence 
	$$
		\dim \mathfrak{isom}(\widetilde{M}^n, \widetilde{g}) = \dim \he(n-2) + \dim \mathfrak{s} = (2n - 3) + \tfrac{(n-3)(n-4)}{2} = \tfrac{n(n+1)}{2} - (2n - 3).
	$$
	This completes the proof.
\end{proof}

\begin{remark}
	We point out that such self-adjoint operators $A \in \End(V)$ whose two eigenspaces yield a decomposition $V = E_1 \oplus E_2$ with $\dim E_2 = 1$ do not occur 
	in the explicit constructions in \cite{5} for those $\Gamma \subset G$ in \cite{5} producing a 
	compact Lorentzian quotient manifold $\widetilde{M}/\Gamma$ in dimensions $n = 3j + 2$ with $j > 1$. 
	Hence it is unclear if such $A \in \End(V)$ can occur in some	isometric identification $(\widetilde{M}/\Gamma, \pi^*g) \simeq (M^n, g)$ discussed in 
	Section \ref{section:Lorentz} for $n > 5$. \par
		More precisely, the examples in \cite{5} contain operators $A \in \End(V)$ with at most three distinct eigenvalues $\lambda_1,\lambda_2,\lambda_3 \in \R$ with
	$\dim E_1 = \dim E_2 = \dim E_3 = \frac{n - 2}{3}$. Hence, choosing $\lambda_1 = \lambda_2 \neq \lambda_3$, we can construct by Corollary \ref{cor:ev} 
	complete, compact Lorentzian manifolds with essentially parallel Weyl tensor whose full isometry group has dimension at least 
	$\frac{n(n+1)}{2} - (2 \cdot (\frac{n - 2}{3})^2 + n)$.
\end{remark}

\begin{proof}
	Choosing $\lambda_1 = \lambda_2 \neq \lambda_3$, the required $F \in \mathfrak{so}(n - 2)$ has to commute with the operator $A$, 
	whence $F(E_1 \oplus E_2) = E_1 \oplus E_2$ and $F(E_3) = E_3$. Thus it holds
	$$
		\dim \mathfrak{isom}(\widetilde{M}^n, \widetilde{g}) = \dim \he(n-2) + \dim \mathfrak{s} = (2n - 3) + 
		\dim \mathfrak{so}(j) + \mathfrak{so}(2j) = (6j + 1) + \tfrac{5j^2 - 3j}{2} = \tfrac{5j^2 + 9j + 2}{2}.
	$$
	Since $\frac{n(n + 1)}{2} - \tfrac{5j^2 + 9j + 2}{2} = 2 \cdot (\frac{n - 2}{3})^2 + n$ this completes the proof.
\end{proof}

	
		

\bibliographystyle{amsalpha}
\bibliography{Bibliography}

\vspace{1cm}

\noindent\textsc{Daniel Schliebner}\newline
Humboldt-Universität zu Berlin, Institut für Mathematik\newline
Rudower Chaussee 25, Room 1.304, D--12489 Berlin.\newline
E-Mail: {\tt schliebn@mathematik.hu-berlin.de}

\end{document}